\documentclass[12pt]{amsart}
\usepackage{amsfonts}
\usepackage{latexsym,amsmath,amsthm,amssymb,amsfonts}

\numberwithin{equation}{section}
\allowdisplaybreaks

\def\endproof{$\hfill\Box$\\}
\def\s{\,\,\,\,}
\def\R{\mathbb{R}}

\def\mv{2.0ex}

\numberwithin{equation}{section}
\newtheorem{theorem}{Theorem}[section]
\newtheorem{lem}[theorem]{Lemma}
\newtheorem{thm}[theorem]{Theorem}
\newtheorem{pro}[theorem]{Proposition}
\newtheorem{cor}[theorem]{Corollary}

\newtheorem{rem}[theorem]{Remark}

\def\dint{\displaystyle{\int}}
\def\lan{\langle}
\def\ran{\rangle}

\newcounter{Cnumber}

\title[ ]
{\bf Blowup Behavior of Harmonic Maps\\ with Finite Index}
\author[ ]
{Yuxiang Li, \quad Lei Liu,\quad
Youde Wang
}

\dedicatory{\sc (In Memory of Prof. Weiyue  Ding)}

\date{}
\begin{document}
\maketitle

\begin{abstract}
In this paper, we study the blow-up phenomena on the $\alpha_k$-harmonic map sequences with bounded uniformly $\alpha_k$-energy, denoted by $\{u_{\alpha_k}: \alpha_k>1 \quad \mbox{and} \quad \alpha_k\searrow 1\}$, from a compact Riemann surface into a compact Riemannian manifold. If the Ricci curvature of the target manifold is of a positive lower bound
and the indices of the $\alpha_k$-harmonic map sequence with respect to the corresponding $\alpha_k$-energy are bounded, then, we can conclude that, if the blow-up phenomena occurs in the convergence of $\{u_{\alpha_k}\}$ as $\alpha_k\searrow 1$, the limiting necks of the convergence of the sequence consist of finite length geodesics, hence the energy identity holds true. For a harmonic map sequence $u_k:(\Sigma,h_k)\rightarrow N$, where the conformal class defined by $h_k$ diverges, we also prove some similar results.
\end{abstract}

\section{Introduction}
Let $(\Sigma,g)$ be a compact Riemann surface and $(N,h)$ be an $n$-dimensional smooth compact Riemannian
manifold which is embedded in $\mathbb{R}^K$ isometrically. Usually, we denote the space of Sobolev maps from $\Sigma$
into $N$ by $W^{k,p}(\Sigma, N)$, which is defined by
$$W^{k,p}(\Sigma, N)=\{u\in W^{k,p}(\Sigma, \mathbb{R}^K): u(x)\in N\,\, \text{for a.e.}\,\,x\in\Sigma\}.$$
For $u\in W^{1,2}(\Sigma,N)$, we define locally the energy
density $e(u)$ of $u$ at $x\in \Sigma$ by
$$e(u)(x)=|\nabla_g u|^2=g^{ij}(x)h_{\alpha\beta}(u(x))
\frac{\partial u^\alpha}{\partial x^i}\frac{\partial
u^\beta}{\partial x^j}.$$
The energy of $u$ on $\Sigma$, denoted by $E(u)$ or $E(u, \Sigma)$, is defined by
$$E(u)=\frac{1}{2}\dint_\Sigma e(u)dV_g,$$
and the critical points of $E$ are called harmonic maps.\medskip

We know that the energy functional $E$ does not satisfy the Palais-Smale condition. In order to overcome this difficulty, Sacks and Uhlenbeck \cite{Sacks-Uhlenbeck1} introduced
the so called $\alpha$-energy $E_\alpha$ of $u: \Sigma\rightarrow N$ as the following
$$E_\alpha(u)=\frac{1}{2}\int_\Sigma\{(1+|\nabla u|^2)^\alpha-1\} dV_g,$$
where $\alpha>1$. The critical points of $E_\alpha$ in $W^{1,2\alpha}(\Sigma,N)$ are called as the $\alpha$-harmonic maps from $\Sigma$ into $N$. It is well-known that this $\alpha$-energy functional $E_\alpha$ satisfies the Palais-Smale condition and therefore there always exists an $\alpha$-harmonic maps in each homotopic class of map from $\Sigma$ into $N$.

The strategy of Sacks and Uhlenbeck is to employ such a sequence of $\alpha_k$-harmonic maps to approximate a harmonic map as $\alpha_k$ tends decreasingly to 1. If the convergence of the sequence of $\alpha_k$-harmonic map is smooth, the limiting map is just a harmonic map from $\Sigma$ into $N$.

The energy of a map $u$ from a closed Riemann surface $\Sigma$ is of conformal invariance, it means that,
if we let $g'=e^{2\varphi}g$ be another conformal metric of $\Sigma$, then
 $$\int_{\Sigma}|\nabla_gu|^2d\mu_g=\int_{\Sigma}|\nabla_{g'}u|^2d\mu_{g'}.$$
Let $\mathcal{C}_g$ denote the conformal class induced by a metric $g$, then, the following definition
$$E(u,\mathcal{C}_g)=\frac{1}{2}\int_\Sigma|\nabla_g u|^2d\mu_g$$
does make sense. Moreover, it is well-known that the critical points of $E(u,\mathcal{C}_g)$ are some branched minimal immersions (see \cite{Sacks-Uhlenbeck1, Sacks-Uhlenbeck2}). Hence, in order to get a branched minimal surface, we also need to study the convergence behavior of a sequence of harmonic maps $u_{k}:(\Sigma,h_k)\rightarrow N$ with bounded uniformly energy $E(u_{\alpha_k}) <C$.\medskip

No doubt, it is very important to study the convergence of a sequence of $\alpha_k$-harmonic maps from a fixed Riemann surface $(\Sigma, g)$ and a sequence of harmonic maps from
$(\Sigma,h_k)$ into $N$, where $h_k$ is the metric with constant curvature. In fact, these problems on the convergence of harmonic map or approximate harmonic map sequences have been studied extensively by many mathematicians. Although these sequences converge smoothly harmonic maps under some suitable geometric and topological conditions, generally one found that the convergence of such two classes of sequences might blow up.\medskip

First, let's recall the convergence behavior of $\alpha_k$-harmonic map sequences. Suppose that $\{u_{\alpha_k}\}$ be a sequence of $\alpha_k$-harmonic maps from $(\Sigma, g)$ with bounded uniformly $\alpha_k$-energy, i.e. $E_{\alpha_k}(u_{\alpha_k})\leq\Theta$. By the theory of Sacks-Uhlenbeck, there exists a subsequence of $\{u_{\alpha_k}\}$, still denoted by $\{u_{\alpha_k}\}$, and a finite set $\mathcal{S}\subset\Sigma$ such that the subsequence converges to a harmonic map $u_0$ in $C^\infty_{loc}(\Sigma\setminus \mathcal{S})$. We know that, at each point $p_i\in\mathcal{S}$, the energy of the subsequence concentrates and the blow-up phenomena occur. Moreover, there exist point sequences $\{x_{i_k}^l\}$ in $\Sigma$ with $\lim\limits_{k\rightarrow+\infty}x_{i_k}^l= p_i$ and scaling constant number sequences $\{\lambda_{i_k}^l\}$ with
 $$\lim\limits_{k\rightarrow+\infty}\lambda_{i_k}^l\rightarrow 0,\s\s l=1,\cdots, n_0,$$
such that
$$u_{\alpha_k}(x_{i_k}^l+\lambda_{i_k}^lx)\rightarrow v^l\s\s \text{in}\s\s C^j_{loc}(\mathbb{R}^2 \setminus \mathcal{A}^i),$$
where all $v^i$ are non-trivial harmonic maps from $S^2$ into $N$, and $\mathcal{A}^i\subset\mathbb{R}^2$ is a finite set.
In order to explore and describe the asymptotic behavior of $\{u_{\alpha_k}\}$ at each blow-up point, the following two problems were raised naturally.\medskip

One is whether or not the energy identity, which states that all the concentrated energy can be accounted for by harmonic bubbles, holds true, i.e.,
$$\lim_{\alpha_k\rightarrow 1}E_{\alpha_k}(u_{\alpha_k}, B^\Sigma_{r_0}(p_i))= E(u_0,B^\Sigma_{r_0}(p_i))+ \sum_{l=1}^{n_0}E(v^l).$$
Here, $B^\Sigma_{r_0}(p_i)$ is a geodesic ball in $\Sigma$ which contains only one blow-up point $p_i$.\medskip

The other is whether or not the limiting necks connecting bubbles are some geodesics in $N$ of finite length?\medskip

For a harmonic map sequence $\{u_k\}$ from $(\Sigma, h_k)$ into $(N, h)$, one also encountered the same problems as above.\medskip

Up to now, for both cases one has made considerably great progress in these two problems \cite{Chen-Tian, Ding, Ding-Tian, Hong-Yin, Lamm, Li-Wang, Li-Wang2, Lin-Wang, Parker, Qing, Qing-Tian, Zhou,Zhu}.

In particular, in \cite{Li-Wang} it is proven that if energy concentration does occur, then a generalized  energy identity holds. Moreover, from the view point of analysis some sufficient and necessary conditions were given such that the energy identity holds true. On the other hand, a relation between the blowup radii and the values of $\alpha$ was discovered to ensure the "no neck property". If necks do occur, however, they must converge to geodesics and a example was given to show that there are even some limiting necks (geodesics) of infinite length.\medskip

Generally, the energy identity does not holds true. For the case of harmonic map sequence $u_k: (\Sigma, h_k)\rightarrow (N,h)$ one has found a counter-example for the energy identity in \cite{Parker}. Very recently, in \cite{Li-Wang2} a counter-example for the energy identity was given for the case of $\alpha_k$-harmonic map sequence.

Furthermore, from the study in \cite{Li-Wang, Li-Wang2} we can see that except for $\alpha$, the topology and geometry of the target manifold $(N, h)$ also play an important role in the convergence of $\alpha$-harmonic map sequence from a compact surface. From the viewpoint of differential geometry, it is therefore natural and interesting to find some reasonable geometric and topological conditions on the domain or target manifold such that the energy identity holds. In particular, a natural question is whether or not we can exploit some geometric and topological conditions to ensure the limiting necks are some geodesics of finite length, which implies that the energy identity holds true? For this goal, in this paper we obtain the following two theorems:

\begin{thm}\label{main1}
Let $(\Sigma,g)$ be a closed Riemann surface and $(N,h)$ be a closed Riemannian manifold with $Ric_N>\lambda>0$. Let $\alpha_k\rightarrow 1$ and $\{u_{\alpha_k}\}$ be a sequence of maps from $(\Sigma, g)\rightarrow (N,h)$ such that each $u_k$ is an $\alpha_k$-harmonic map, the indices and energy satisfy respectively
$$\mbox{Index}(E_{\alpha_k}(u_{\alpha_k}))<C,\s \s \s E_{\alpha_k}
(u_{\alpha_k})<C.$$
If $\{u_{\alpha_k}\}$ blows up,  then the limiting necks consist of some finite length geodesics.
\end{thm}

\begin{thm}\label{main2}
Let $\Sigma$ be a closed Riemann surface with genus $g(\Sigma)\geq 1$. In the case $g(\Sigma)\geq 2$, $\Sigma$ is equipped a sequence of smooth metrics $h_k$ with curvature $-1$. In the case $g(\Sigma)=1$,  $\Sigma$ is equipped a sequence of smooth metrics $h_k$ with curvature $0$ and the area $A(\Sigma, h_k) = 1$. Let $(N, h)$ be a Riemannian manifold with the  Ricci curvature $\mbox{Ric}_N >\lambda> 0$. Suppose that $(\Sigma, h_k)$ diverges in the moduli space and $\{u_k\}$ is a harmonic map sequence from $(\Sigma, h_k)$ into $(N, h)$ with bounded index and energy. If the set of the limiting necks of $u_k$ is not empty, then it consists of finite length geodesics.
\end{thm}

\begin{rem}
By the results in \cite{Chen-Tian} or \cite{Li-Wang}, the fact the limiting necks are of finite length implies that the energy identity is true. We should also mention that, when each $u_{\alpha_k}$ in $\{u_{\alpha_k}\}$ is the minimizer of the corresponding $E_{\alpha_k}$ in a fixed homotopy class, Chen and Tian \cite{Chen-Tian} have proved that the necks are just some geodesics of finite length in $N$.
\end{rem}

\begin{rem}
The curvature condition in Theorem \ref{main1} and \ref{main2} is used to ensure that any geodesic of infinite length lying on $N$ is not stable. In fact, we will prove in this paper that, if the necks contain a unstable geodesic of infinite length, then the indices of the harmonic (or $\alpha$-harmonic) map sequence can not be bounded from the above.
\end{rem}\medskip

\section{The Proofs of Theorem \ref{main1}}

Our strategy is to show that the indices of the sequence $\{u_{\alpha_k}\}$ in Theorem 1.1 are not bounded if there exists a infinite length geodesic in the set of the limiting necks $\{u_{\alpha_k}\}$. For this goal, first we need to recall the definition of the index of a $\alpha$-harmonic map and the second variational formula of $\alpha$-energy functional.

\subsection{The index of a $\alpha$-harmonic map}

Let $u:(\Sigma,g)\rightarrow (N,h)$ be an $\alpha$-harmonic map. $L=u^{-1}(TN)$ is a smooth pull-back bundle over $\Sigma$. Let $V$ be a section of $L$ and $$u_t(x)=exp_{u(x)}(tV).$$ Obviously, $u_0=u$. Then, the formula of the second variation of $E_\alpha$ reads
\begin{eqnarray}\label{second.Ealpha}
\delta^2E_\alpha(u)(V,V)&=&\frac{d^2}{dt^2}E_\alpha(u_t)|_{t=0}\nonumber\\[\mv]
&=&2\alpha\dint_{\Sigma}(1+|du|^2)^{(\alpha-1)}\left(\langle \nabla V,\nabla V \rangle-R(V,\nabla u,\nabla u,V)\right)d\mu\\[\mv]
&&+4\alpha(\alpha-1)\dint_\Sigma(1+|du|^2)^{\alpha-2}
\lan du,\nabla V\ran^2d\mu.\nonumber
\end{eqnarray}
For more details we refer to \cite{M-M}.

Let $\Gamma(L)$ denotes the linear space of the smooth sections of $L$. Then, the index of $u$ is the maximal dimension of the linear subspaces of $\Gamma(L)$ on which the \eqref{second.Ealpha} is definite negatively.

\subsection{The index of the necks}
We have known the limiting necks of $\{u_{\alpha_k}\}$ are some geodesics in $N$, a natural question is there exists some relations between the indices of these geodesics and the indices of the necks of $\{u_{\alpha_k}\}$. In this subsection we need to analyse the asymptotic behavior of the necks of  $\{u_{\alpha_k}\}$ and try to establish the desired relations.

Let $\alpha_k\rightarrow 1$ and each $u_{\alpha_k}$ of the map sequence $\{u_{\alpha_k}: k=1,2,\cdots\}$ be a $\alpha_k$ harmonic map from $(\Sigma,g)$ into $(N,h)$. For convenience we always embed $(N,h)$ into $\R^K$ isometrically
and set $u_k=u_{\alpha_k}$. Assume that $\{u_k\}$ blows up only at a point $p\in\Sigma$. Then, for any $\epsilon$, we have
 $$\lim_{k\rightarrow+\infty}\|\nabla u_{k}\|_{C^0(B_\epsilon(p))}=+\infty.$$
Choose an isothermal coordinate chart $(D;x^1,x^2)$ centered at $p$,
such that
$$g=e^{2\varphi}(dx^1\otimes dx^1+dx^2\otimes dx^2),\s \mbox{ and }\s\varphi(0)=0.$$
For simplicity, we assume $u_{k}$ has only one blowup
point in $D$. Put
$$r_k=\frac{1}{\|\nabla u_{k}\|_{C^0(D_\frac{1}{2})}},\s
\mbox{ and }\s
|\nabla u_{k}(x_k)|=\|\nabla u_{k}\|_{C^0(D_\frac{1}{2})}.$$
Then, we have that $x_k\rightarrow 0$, $r_k\rightarrow 0$ and there exists a bubble $v$, which can be considered as a harmonic map from $S^2$ into $N$, such that $u_k(x_k+r_kx)$ converges to $v$. Without loss of generality, we may assume $x_k=0$.\medskip

By the arguments in \cite{Li-Wang}, we only need to prove Theorem \ref{main1} for the case there exists one bubble in the convergence of $\{u_k\}$. So, we always assume that only one bubble appears in the convergence of $\{u_{k}\}$ in this section. \medskip

Now, we consider the case that the limiting necks contain a geodesic of infinite length. In fact, the present paper is a follow up of the papers \cite{Li-Wang} and \cite{Chen-Li-Wang}, first of all, we need to recall some results proved in \cite{Li-Wang}.

\begin{lem} Let $\alpha_k\rightarrow 1$ and $\{u_{k}\}$ be a map sequence such that each $u_{k}$ is an $\alpha_k$-harmonic map from $(\Sigma,g)$ into $(N,h)$. If there is a positive constant $\Theta$ such that $E_{\alpha_k}(u_{k})<\Theta$ for any $\alpha_k$, then, there exists a positive constant $C$ such that, neglecting a subsequence, there holds
$$\|\nabla u_{k}\|^{\alpha_k-1}_{C^0(\Sigma)}<C.$$
\end{lem}
For the proof of this lemma and more details we refer to Remark 1.2 in \cite{Li-Wang}. Moreover, for the convergence radii and $\alpha_k$ we have following relations:

\begin{lem} Let $\{u_{k}\}$ satisfy the same conditions as in Lemma 2.1. If there exists only one bubble in the convergence of $\{u_{k}\}$ and the limiting neck is of infinite length, then, the following hold true
$$0<-(\alpha_k-1)\log r_k<C,
\s\mbox{and}\s\s \sqrt{\alpha_k-1}\log r_k\rightarrow-\infty.$$
Here, $r_k$ is defined as before.
\end{lem}

\proof From Remark 1.2 in \cite{Li-Wang}, we have $\mu=\liminf_{\alpha_k\rightarrow 1}r_k^{2-2\alpha_k}\in [1,\,\mu_{\max}]$ where $\mu_{\max}\geq 1$ is a positive constant. Therefore, it follows that there holds
$$0<-(\alpha_k-1)\log r_k<C.$$

Since the limiting neck is of infinite length, from Theorem 1.3 in \cite{Li-Wang} we known that
$$\nu=\liminf_{\alpha_k\rightarrow 1}r_k^{-\sqrt{\alpha_k-1}}\rightarrow\infty.$$
It follows that  $$\sqrt{\alpha_k-1}\log r_k\rightarrow-\infty.$$
Thus we complete the proof.\endproof

As a direct corollary of the Proposition 4.3 in \cite{Li-Wang}, we have
\begin{lem}\label{key} Let $\alpha_k\rightarrow 1$ and $\{u_{k}\}$ be a map sequence such that each $u_{k}$ is an $\alpha_k$-harmonic map from $(\Sigma,g)$ into $(N,h)\subset\mathbb{R}^K$. Suppose that there is a positive constant $\Theta$ such that $E_{\alpha_k}(u_{k})<\Theta$ for any $\alpha_k$ and there exists only one bubble in the convergence of $\{u_{k}\}$. Then, for any $t_k\rightarrow t\in(0, 1)$, there exist a vector $\xi\in\mathbb{R}^K$ and a subsequence of  $\{u_k\}$ such that
\begin{equation}
\frac{1}{\sqrt{\alpha_k-1}}\frac{\partial u_{k}}{\partial \theta}(r_k^{t_k}e^{\sqrt{-1}\theta})\to 0
\end{equation}
and
\begin{equation}
\frac{r_k^{t_k}}{\sqrt{\alpha_k-1}}\frac{\partial u_{k}}{\partial r}(r_k^{t_k}e^{\sqrt{-1}\theta})\to \xi
\end{equation}
as $k\to \infty$. Moreover,
$$|\xi|=\mu^{1-t}\sqrt{\frac{E(v)}{\pi}},$$
where $\mu$ is defined by
$$\mu=\lim_{k\rightarrow+\infty}r_k^{2-2\alpha_k}.$$
\end{lem}

Now we define the approximate curve of $u_k$, denoted by $u_k^*(r)$, by
 $$u_k^*(r)=\frac{1}{2\pi}\int_0^{2\pi}u_k(re^{\sqrt{-1}\theta})d\theta.$$
Since the target manifold $(N, h)$ is embedded in $\mathbb{R}^K$, $u_k^*(r)$ is a space curve of $\mathbb{R}^K$ and we denote the arc-length parametrization of $u_k^*$ by $s$ such that $s(r_k^{t_1})=0$. Then

\begin{lem}
Let $\{u_{k}\}$ satisfy the same conditions as in Lemma 2.3. Suppose that the limiting neck of $\{u_{k}\}$ is a geodesic of infinite length. Then, there exists a subsequence of $\{u^*_k(s)\}$ which converges smoothly on $[0,a]$ to a geodesic $\gamma$ for any fixed $a>0$.
\end{lem}

Without loss of generality, from now on, we assume that
$u_k^*(s)$ converges smoothly to $\gamma$ on any $[0,a]$.
As a corollary, we have

\begin{cor}\label{convergence.of.uk}
Let $\{u_{k}\}$ satisfy the same conditions as in Lemma 2.3. Suppose that the limiting neck of $\{u_{k}\}$ is a geodesic of infinite length. Then, for any given $a>0$ and any fixed $\theta$, $u_k(se^{\sqrt{-1}\theta})$ converges to $\gamma$
in $C^1[0,a]$. Moreover, we have
\begin{equation}\label{arclength}
\left\|\frac{r(s)}{\sqrt{\alpha_k-1}}\left|\frac{\partial s}{\partial r}\right|-\mu^{1-t_1}\sqrt{\frac{E(v)}{\pi}}\right\|_{C^0([0,a])}\rightarrow 0.
\end{equation}
\end{cor}

\proof Let $$s(r_k^{t_k^a})=a.$$ By Lemma \ref{key}, we have
\begin{equation}\label{theta}
a=\int_{r_k^{t_k^a}}^{r_k^{t_1}}\left|\frac{d{u}_k^*(r)}{dr}\right|dr\geq
C\int_{r_k^{t_k^a}}^{r_k^{t_1}}\frac{\sqrt{\alpha_k-1}}{r}dr=
-C(t_k^a-t_1)\sqrt{\alpha_k-1}\log r_k.
\end{equation}
On the other hand, Lemma 2.2 (see Theorem 1.3 in \cite{Li-Wang}) tells us
 $$\sqrt{\alpha_k-1}\log r_k\rightarrow -\infty,$$
since the limiting neck (geodesic) is of infinite length. Hence, from (\ref{theta}) and the above fact, we have
\begin{equation}\label{ta}
t_k^a-t_1\rightarrow 0\s\s\mbox{ as }\s\s k\rightarrow+\infty.
\end{equation}

We assume that $u_k(se^{\sqrt{-1}\theta})$ does not converge to $\gamma$ in $C^1[0,a]$. Then there exists $s_{k_i}\in [0,a]$, such that
$$\sup_{\theta}\left|\frac{\partial u_k}{\partial s}(s_{k_i}e^{\sqrt{-1}\theta})-\frac{d{u}_k^*}{ds}(s_{k_i})\right|>\epsilon>0.$$
Let $s_{k_i}=r_k^{t_{k_i}}$. Obviously, $t_{k_i}\in[t_1,\,t_{k_i}^a]$. Thus $t_{k_i}\rightarrow t_1$. By Lemma \ref{key}, after passing to a subsequence, we have
$$\lim_{k\rightarrow\infty}\frac{r_k^{t_{k_i}}}{\sqrt{\alpha_k-1}}\left|\frac{\partial u_k}{\partial r}(r_k^{t_{k_i}}e^{\sqrt{-1}\theta})-\frac{d{u_k}^*}{dr}(r_k^{t_{k_i}})
\right|\rightarrow 0.$$
Therefore, noting
$$\left|\frac{\partial u_k(s_{k_i}e^{\sqrt{-1}\theta})}{\partial s}-\frac{d{u}_k^*(s_{k_i})}{ds}\right|=
\left|\frac{dr}{ds}\right|\cdot\left|\frac{\partial u_k}{\partial r}-\frac{d{u_k}^*(r)}{dr}\right|_{r=r_k^{t_i}}$$
and $$\left|\frac{dr}{ds}\right|_{r=r_k^{t_{k_i}}}\leq \frac{Cr_k^{t_{k_i}}}{\sqrt{\alpha_k-1}},$$
we have
$$\left|\frac{\partial u_k(s_{k_i}e^{\sqrt{-1}\theta})}{\partial s}-\frac{d{u}_k^*(s_{k_i})}{ds}\right| \leq \frac{Cr_k^{t_i}}{\sqrt{\alpha_k-1}}\left|\frac{\partial u_k}{\partial r}(re^{\sqrt{-1}\theta})-\frac{d{u_k}^*}{dr}(r)\right|_{r=r_k^{t_i}}\rightarrow 0.$$
Thus, we get a contradiction. Hence, it follows
$$\left\|\frac{\partial u_k(se^{\sqrt{-1}\theta})}{\partial s}-\frac{d{u}_k^*(s)}{ds}\right\|_{C^0[0,\,a]}\rightarrow 0.$$
From the arguments in the above and \cite{Li-Wang} we conclude that for any fixed $\theta$
$$\|u_k(se^{\sqrt{-1}\theta})-u_k^*(s)\|_{C^1[0,\,a]}\rightarrow 0.$$
By the same way, we can prove \eqref{arclength}.
~\endproof

\begin{lem}\label{LW} Suppsoe that $\{u_{k}\}$ satisfies the same conditions as in Lemma 2.3. Then, for any fixed $R>0$ and $0<t_1<t_2<1$, we have
 $$\lim_{k\rightarrow+\infty}\sup_{t\in[t_1,t_2]}\frac{1}{\alpha_k-1}\int_{D_{Rr_k^{t}}\setminus D_{\frac{1}{R}r_k^{t}}} |u_{k,\theta}|^2dx=0,$$
where $$u_{k,\theta}=r^{-1}\frac{\partial u_k}{\partial\theta}.$$
\end{lem}

\proof Assume this is not true. After passing to a subsequence,
we can find $t_k\rightarrow
t\in [t_1,t_2]$, such that
$$
\frac{1}{\alpha_k-1}\int_{D_{Rr_k^{t_k}}\setminus D_{\frac{1}{R}r_k^{t_k}}}
|u_{k,\theta}|^2dx\geq\epsilon.
$$
However, by Proposition 4.2 in \cite{Li-Wang},
$$
\lim_{k\rightarrow+\infty}\frac{1}{\alpha_k-1}\int_{D_{Rr_k^{t_k}}\setminus D_{\frac{1}{R}r_k^{t_k}}}
|u_{k,\theta}|^2dx=0.
$$
This is a contradiction, Thus we complete the proof of the lemma.
\endproof

Now, let's recall the definition of stability of a geodesic on a Riemannian manifold $(N,h)$. A geodesic $\gamma$ is called unstable if and only if the second variation formula of its length satisfies
$$I_\gamma(V_0,V_0)=\int_0^a(\langle \nabla_{\dot{\gamma}}V_0,\nabla_{\dot{\gamma}}V_0 \rangle-R(V_0,\dot{\gamma},\dot{\gamma},V_0))ds<-\delta<0.$$
Here $R$ is the curvature operator of $N$. We have

\begin{lem}
Suppose that $\{u_{k}\}$ satisfies the same conditions as in Lemma 2.3. If the limiting neck of $\{u_{k}\}$ is a unstable geodesic which is parameterized on $[0,\, a]$ by arc length, then, for sufficiently large $k$, there exists a section $V_k$ of $u_{k}^{-1}(TN)$, which is supported in $D_{r_k^{t_1}}\setminus D_{r_k^{t_k^a}}(x_k)$, such that
$$\delta^2E_{\alpha_k}(V_k,V_k)<0.$$
\end{lem}

\proof Since the limiting neck of $\{u_{k}\}$, denoted by $\gamma:[0,\, a]\rightarrow N$, is not a stable geodesic, there exists a vector field $V_0$ on $\gamma$
 with $V_0|_{\gamma(0)}=0$ and $V_0|_{\gamma(a)}=0$ such that
 $$I_\gamma(V_0, V_0)<0.$$

Let $P$ be projection from $T\R^K$ to $TN$. We define
$$V_k(t(s)e^{\sqrt{-1}\theta}+x_k)=P_{u_k(se^{\sqrt{-1}\theta})}(V_0(s)),$$ where $s$ is the arc-length parametrization of $u_{k}^*(t)$ with $s(r_k^{t_1})=0$.
Then, $V_k$ is smooth section of $u_{k}^{-1}(TN)$ which
is supported in $D_{r_k^{t_1}}\setminus D_{r_k^{t_k^a}}(x_k)$.
By Corollary \ref{convergence.of.uk}, for any fixed $\theta$, we have that $V_k(u_k(se^{\sqrt{-1}\theta}))$ converges to
$V_0(\gamma(s))$ in $C^1[0,a]$. Then
\begin{eqnarray}\label{second.Ealpha1}
\delta^2E_{\alpha_k}(V_k,V_k)&=&
2\alpha_k\dint_{D_{r_k^{t_1}}\setminus
D_{r_k^{t_k^a}}}(1+|du_k|^2)^{(\alpha_k-1)}\left(\langle \nabla V_k,\nabla V_k\rangle-R(V_k,\nabla u_k,\nabla u_k,V_k)\right)dx\nonumber\\[\mv]
&&+4\alpha_k(\alpha_k-1)\dint_{D_{r_k^{t_1}}\setminus
D_{r_k^{t_k^a}}}(1+|du_k|^2)^{\alpha_k-2}\lan du_k,\nabla V_k\ran^2dx.
\end{eqnarray}
Next, we will show that
\begin{equation}\label{I.I2}
\lim_{k\rightarrow+\infty}
\frac{1}{\sqrt{\alpha_k-1}}\delta^2 E_{\alpha_k}(V_k,V_k)
=4\pi\mu \sqrt{\frac{E(v)}{\pi}}I_\gamma(V_0,V_0).
\end{equation}

We compute
\begin{eqnarray*}
& &\delta^2E_{\alpha_k}(V_k,V_k)\\
&=&2\alpha_k\int_0^{2\pi}\int^{r_k^{t_1}}_{r_k^{t_k^a}}
(1+|du_k|^2)^{(\alpha_k-1)}(\langle \nabla V_k,\nabla V_k\rangle-R(V_k,\nabla u_k,\nabla u_k,V_k)) rdrd\theta\\
&&+4\alpha_k(\alpha_k-1)\dint_{D_{r_k^{t_1}}\setminus
D_{r_k^{t_k^a}}}(1+|du_k|^2)^{\alpha_k-2}\lan du_k,\nabla V_k\ran^2dx\\
&=&
2\alpha_k\int_0^{2\pi}\int^{r_k^{t_1}}_{r_k^{t_k^a}}
(1+|du_k|^2)^{(\alpha_k-1)}(\langle \nabla_{\frac{\partial u_k}{\partial r}} V_k,\nabla_{\frac{\partial u_k}{\partial r}} V_k\rangle
-R(V_k,\frac{\partial u_k}{\partial r},\frac{\partial u_k}{\partial r},V_k)) rdrd\theta\\
&&+2\alpha_k\int_0^{2\pi}\int^{r_k^{t_1}}_{r_k^{t_k^a}}(1+|du_k|^2)^{(\alpha_k-1)}
(\langle \nabla_{ u_{k,\theta}} V_k,\nabla_{u_{k,\theta}} V_k \rangle-R(V_k,u_{k,\theta},u_{k,\theta},V_k)) rdrd\theta\\
&&+4\alpha_k(\alpha_k-1)\dint_{D_{r_k^{t_1}}\setminus
D_{r_k^{t_k^a}}}(1+|du_k|^2)^{\alpha_k-2}\lan du_k,\nabla V_k\ran^2 dx\\
&=&
2\alpha_k\mathbf{I}+2\alpha_k\mathbf{II}+4\alpha_k\mathbf{III}.
\end{eqnarray*}

Firstly, we calculate $\mathbf{I}$:
\begin{eqnarray*}
& &\frac{\mathbf{I}}{\sqrt{\alpha_k-1}}\\
&=&\int_0^{2\pi}\int_{0}^{a}(1+|du_k|^2)^{(\alpha_k-1)}\left(\langle \nabla_{\frac{\partial u_k}{\partial s}} V_k,\nabla_{\frac{\partial u_k}{\partial s}} V_k\rangle
-R(V_k,\frac{\partial u_k}{\partial s},\frac{\partial u_k}{\partial s},V_k)\right)\frac{\left|\frac{\partial s}{\partial r}\right|}{\sqrt{\alpha_k-1}} rdsd\theta.
\end{eqnarray*}
By Lemma 2.3 we can see easily that
$$\left(\left|\frac{r_k^{t(s)}}{\sqrt{\alpha_k-1}}du_k\right|^2
\frac{\alpha_k-1}{r_k^{2t(s)}}\right)^{(\alpha_k-1)}\longrightarrow \mu^{t_1}.$$
It follows from the fact $\mu\geq 1$ (see \cite{Li-Wang}) and the above
$$(1+|du_k|^2)^{(\alpha_k-1)}=\left(1+\left|\frac{r_k^{t(s)}}{\sqrt{\alpha_k-1}}du_k\right|^2
\frac{\alpha_k-1}{r_k^{2t(s)}}\right)^{(\alpha_k-1)}\longrightarrow\mu^{t_1}.$$
Hence, we infer from the above and Corollary \ref{convergence.of.uk}
$$\lim_{k\rightarrow+\infty}\frac{\mathbf{I}}{\sqrt{\alpha_k-1}}
=2\pi\mu\sqrt{\frac{E(v)}{\pi}}I_\gamma(V_0,V_0).$$

\medskip
Next, we calculate the term $\mathbf{II}$. By the definition we have
$$\nabla_{\frac{\partial u_k}{\partial\theta}}V_k = P_{u_k(se^{\sqrt{-1}\theta})}\left(\frac{\partial V_k}
{\partial \theta}\right)=P_{u_k(se^{\sqrt{-1}\theta})}\left(\frac{\partial }{\partial\theta}(P_{u_k(se^{\sqrt{-1}\theta})})(V_0)\right),$$
where $\frac{\partial V_k}{\partial\theta}$ is the derivative
in $\R^n$.
This leads to
$$|\nabla_{\frac{\partial u_k}{\partial\theta}}V_k|\leq C(a)
\left|\frac{\partial u_k}{\partial\theta}\right|.$$
Hence, we have
\begin{eqnarray*}
\frac{\mathbf{II}}{\sqrt{\alpha_k-1}}
&=&
\int_0^{2\pi}\int^{r_k^{t_1}}_{r_k^{t_k^a}}\frac{(1+|du_k|^2)^{(\alpha_k-1)}}{\sqrt{\alpha_k-1}}
\left(\langle \nabla_{ u_{k,\theta}} V_k,\nabla_{u_{k,\theta}} V_k \rangle-R(V_k,u_{k,\theta},u_{k,\theta},V_k)\right)rdrd\theta\\
&\leq&\frac{C}{\sqrt{\alpha_k-1}}\int_0^{2\pi}\int^{r_k^{t_1}}
_{r_k^{t_k^a}}|u_{k,\theta}|^2rdrd\theta.
\end{eqnarray*}
For a given $R>0$, set
$$m_k=\left[\frac{\log r_k^{t_1-t_k^a}}{\log R}\right]+1.$$
It is easy to see that
$$D_{r_k^{t_1}}\setminus D_{r_k^{t_k^a}}\subset\cup_{i=1}^{m_k}(D_{R^ir_k^{t_k^a}}\setminus D_{R^{i-1}r_k^{t_k^a}}).$$
By \eqref{theta}, we have
$$\sqrt{\alpha_k-1}m_k\leq C(R).$$
Then
\begin{eqnarray*}
\frac{\mathbf{II}}{\sqrt{\alpha_k-1}}
&\leq&
\frac{C}{\sqrt{\alpha_k-1}}\int_0^{2\pi}\int^{r_k^{t_1}}
_{r_k^{t_k^a}}
|u_{k,\theta}|^2rdrd\theta\\
&\leq&
\frac{C}{\sqrt{\alpha_k-1}}\int_{\cup_{i=1}^{m_k}(D_{R^ir_k^{t_k^a}}\setminus D_{R^{i-1}r_k^{t_k^a}})}
|u_{k,\theta}|^2dx\\
&\leq&
\frac{Cm_k\sqrt{\alpha_k-1}}{\alpha_k-1}\frac{1}{m_k}\int_{\cup_{i=1}^{m_k}(D_{R^ir_k^{t_k^a}}\setminus D_{R^{i-1}r_k^{t_k^a}})}
|u_{k,\theta}|^2dx\\
&\leq&
\frac{C(R)}{m_k}\left(\frac{1}{\alpha_k-1}\int_{\cup_{i=1}^{m_k}(D_{R^ir_k^{t_k^a}}\setminus D_{R^{i-1}r_k^{t_k^a}})}
|u_{k,\theta}|^2dx\right).
\end{eqnarray*}
It follows from Lemma \ref{LW} and the above inequality that there holds
$$
\lim_{k\rightarrow\infty}\frac{1}{\sqrt{\alpha_k-1}}\mathbf{II}=0.
$$

Lastly, we consider the term $\mathbf{III}$. It is easy to check that
$$|\lan du_k,V_k\ran|\leq C|du_k|^2.$$
So, there exists a constant $C$ such that
$$(1+|du_k|^2)^{\alpha_k-2}\lan du_k,\nabla V_k\ran^2\leq
(1+|du_k|^2)^{\alpha_k-1}<C.$$
This leads to
\begin{eqnarray*}
\frac{\mathbf{III}}{\sqrt{\alpha_k-1}}\leq C\sqrt{\alpha_k-1}
\int_{D_{r_k^{t_1}}}(1+|du_k|^2)^{\alpha_k-1}dx\rightarrow 0.
\end{eqnarray*}
Thus, we obtain the desired estimate and finish the proof.
\endproof

Since that $(N, h)$ is a complete Riemannian manifold with  $\mbox{Ric}_N\geq \lambda>0$, then, the well-known Myers theorem tells us that the diameter of $(N, h)$ satisfies $$\mbox{diam}(N, h)\leq \frac{\pi}{\sqrt{\lambda(n-1)^{-1}}},$$
and any geodesic $\gamma$ lying on $(N, h)$ is unstable if its length $l(\gamma)$ satisfies
 $$l(\gamma)\geq\frac{\pi}{\sqrt{\lambda(n-1)^{-1}}}\equiv l_N.$$
Hence, for any given positive number $a$ such that $a \geq l_N+2\epsilon$ and any geodesic $\gamma$ lying on $(N, h)$ which is parameterized by arc-length in $[0,a]$, there always exists a vector field $V_0(s)$, which is smooth on $\gamma$, and 0 on $\gamma|_{[0,\,\epsilon]}$ and $\gamma|_{[a-\epsilon,\,a]}$,
such that the second variation of the length of $\gamma$ satisfies
\begin{equation}\label{d2.geodesic}
I_\gamma(V_0,V_0)=\int_0^a(\langle \nabla_{\dot{\gamma}}V_0,\nabla_{\dot{\gamma}}V_0 \rangle-R(V_0,\dot{\gamma},\dot{\gamma},V_0))ds<-\delta<0.
\end{equation}

\begin{lem}
Let $(N, h)$ be a closed Riemannian manifold with $\mbox{Ric}(N)\geq\lambda>0$. Suppose that $\{u_{k}\}$ satisfies the same conditions as in Lemma 2.3. If the limiting neck of $\{u_{k}\}$ is a geodesic of infinite length, then the indices of $\{u_{k}\}$ with respect to the corresponding $E_{\alpha_k}$ can not be bounded from above.
\end{lem}

\proof
Since the limiting neck of $\{u_{k}\}$ is a geodesic of infinite length, then, for given $t_1$, the above arguments in Lemma 2.7 tell us that we can always choose a suitable positive constant $\epsilon_1$ such that, as $k$ is large enough, the arc length $a$ of $u^*_k(s)$ on $D_{r_k^{t_1}}\setminus D_{r_k^{t_1+\epsilon}}(x_k)$ satisfies
$$a > l_N=\frac{\pi}{\sqrt{\lambda(n-1)^{-1}}}.$$
Therefore, there exists a section $V^1_k$ of $u_k^{-1}(TN)$, which is $0$ outside $D_{r_k^{t_1}}\setminus D_{r_k^{t_1+\epsilon}}(x_k)$, satisfying
$$\delta^2E_{\alpha_k}(V_k^1,V_k^1)<0.$$
By the same method, for $t_2=t_1+2\epsilon_1$, we can also pick $\epsilon_2>0$ and construct a section $V_k^2$, which is $0$ outside
$D_{r_k^{t_2}}\setminus D_{r_k^{t_2+\epsilon_2}}(x_k)$, such that
$$\delta^2E_{\alpha_k}(V_k^2,V_k^2)<0.$$
Since the limiting neck is a geodesic of infinite length,  then, when $k$ is sufficiently large, there exists $i_k$ with $i_k\rightarrow\infty$ such that we can construct by the same way as above a series of sections $\{V_k^3, V_k^4, \cdots, V_k^{i_k}\}$ satisfying that for any $1\leq i\leq i_k$ there holds true
$$\delta^2E_{\alpha_k}(V_k^i,V_k^i)<0.$$
Obviously, $V_k^1$, $V_k^2$, $\cdots$, $V_k^{i_k}$ are linearly independent. This means that $$\mbox{Index}(E_{\alpha_k}(u_{k}))\geq i_k. $$
Therefore, we get
$$\mbox{Index}(E_{\alpha_k}(u_{k}))\rightarrow+\infty, \s\s\s\mbox{ as }\s k\rightarrow+\infty.$$
Thus, we complete the proof of the lemma.
\endproof
\medskip

\noindent{\bf The proof of Theorem 1.1}: Obviously, Theorem \ref{main1} is just a direct corollary of the above lemma.
\medskip

\section{The Proofs of Theorem \ref{main2}}

From the arguments and the appendix in \cite{Chen-Li-Wang} we know that one only need to consider the convergence behavior of harmonic map sequences from two dimensional flat cylinders, although the original harmonic map sequence is from a sequence of hyperbolic or flat closed Riemann surfaces respectively. First, we recall some fundamental notions such as the index of a harmonic map with respect to the energy functional. \medskip

Let $T_k\rightarrow\infty$ be a series of positive numbers and $u:(-T_k,T_k)\times S^1\rightarrow (N,h)$ be a harmonic map. $L=u^{-1}(TN)$ is the pull-back bundle over $(-T_k,T_k)\times S^1$. Let $V$ be a section of $L$ which is 0 near $\{\pm T_k\}\times S^1$ and
 $$u_\tau(x)=exp_{u(x)}(\tau V).$$
It is well-known that the second variational formula of energy functional $E$ is the following:

\begin{eqnarray*}
\delta^2E(u)(V,V)
&=&2\dint_{\Sigma}\left(\langle \nabla V,\nabla V \rangle-R(V,\nabla u,\nabla u,V)\right)dtd\theta.
\end{eqnarray*}
Let $\Gamma(L)$ denote the linear space of the smooth sections of $L$. Then, the index of $u$ is just the maximal dimension of a linear subspace of $\Gamma(L)$ on which the above is definite negatively.\medskip

Let $u_k$ be an harmonic map from $(-T_k,T_k)\times S^1$ into $(N,h)$. We assume that, for any $t_k\in (-T_k,T_k)$,
$$|\nabla u_k(\theta,t_k)|\rightarrow 0,\s\s\s \mbox{as}\s k\rightarrow\infty.$$
Moreover, we assume that $u_k((-T_k,T_k)\times S^1)$ converges to an infinite length geodesic.\medskip

By the arguments in \cite{Chen-Li-Wang}, we can see easily that Theorem \ref{main2} in this paper can be deduced from the following lemma:

\begin{pro}\label{sy}
Let $\{u_k:(-T_k,T_k)\times S^1\rightarrow N, k=1, 2, \cdots\}$ be a harmonic map sequence such that for any $t_k\in (-T_k,T_k)$, there holds true $|\nabla u_k(\theta,t_k)|\rightarrow 0$. If $\mbox{Ric}_N\geq\lambda>0$ and $u_k((-T_k,T_k)\times S^1)$ converges to an infinite length geodesic, then the index of $u_k$ tends to infinity.
\end{pro}

In order to prove Proposition \ref{sy}, we need to recall some known results which were established in \cite{Chen-Li-Wang}. We first recall a useful observation in \cite{Zhu}.

\begin{lem}\label{zhu1}Let $u$ be a harmonic map from $(-T,\, T)\times S^1
\rightarrow N$. Then, the following function defined by
$$\beta(u)=\int_{\{t\}\times S^{1}}(|u_{t}|^{2}-|u_{\theta}|^{2}-2iu_{t}\cdot u_{\theta})d\theta$$
is independent of $t\in(-T,\, T)$.
\end{lem}

Next, we recall some known results proved in \cite{Chen-Li-Wang}, which are used in the following arguments.

\begin{lem} Let $\{u_k:(-T_k,\,T_k)\times S^1\rightarrow N, \,\, k=1, 2, \cdots\}$ be a sequence of harmonic maps such that for any $t_k\in (-T_k,\,T_k)$, there holds true $|\nabla u_k(\theta,t_k)|\rightarrow 0$. Assume that $u_k((-T_k,\,T_k)\times S^1)$ converges to an infinite length geodesic. Then, as $k\rightarrow 0$, we have
$$\lim_{k\rightarrow\infty}\sqrt{|\mbox{Re} \ \beta(u_{k})|}T_{k}
=\infty.$$
\end{lem}

\begin{lem}\label{key2}
Let $\{u_k\}$ satisfy the same conditions as in Lemma 3.3. Then, for any $\lambda<1$ and $t_k\in [-\lambda T_k,\,\lambda T_k]$, there exists a vector $\xi\in\mathbb{R}^K$ and a subsequence of
$$\left\{\frac{1}{\sqrt{| \mbox{Re}(\beta(u_{k}))|}}\frac{\partial
u_{k}}{\partial t}(t_k, \theta):\,\, k=1, 2, 3, \cdots\right\}$$ such that the subsequence converges to $\xi$. Moreover, we have
$$|\xi|=\frac{1}{\sqrt{2\pi}}.$$
\end{lem}

By Lemma 2.6 in \cite{Chen-Li-Wang}, we also have
\begin{lem}\label{theta.energy2}
Let $\{u_k\}$ satisfy the same conditions as in Lemma 3.3. Then, for any fixed $0<\lambda<1$ and  $T>0$, we have
$$\lim_{k\rightarrow\infty}\sup_{t\in [-\lambda T_k,\,\,\lambda T_k]}\frac{1}{|\mbox{Re} \ \beta(u_{k})|}\int_{[t-T,\,\, t+T]\times S^1}\left|\frac{\partial u_k}{\partial\theta}\right|^{2}dtd\theta=0.$$
\end{lem}

\medskip
As in \cite{Chen-Li-Wang}, we introduce the following sequence of curves in $\mathbb{R}^K$ defined by
$$u_k^*(t)=\frac{1}{2\pi}\int_0^{2\pi}u_k(t,\theta)d\theta.$$
Obviously, these curves are smooth. Now, for each $k$, let $s$ be the arc-length parametrization of the the curve $u_k^*(t)$ with $s(0)=0$.
By the arguments in \cite{Chen-Li-Wang}, we have

\begin{lem}
Under the same conditions as Lemma 3.3, we have that $u^*_k(s)$
converges smoothly on $[0,a]$ to a geodesic $\gamma$ on $N$ for any fixed $a>0$.
\end{lem}

From now on, we assume $u_k^*(s)$ converges to $\gamma$
on $[0,a]$ for any fixed $a>0$.
Set $s(t_k^a)=a$. Similar to Corollary \ref{convergence.of.uk}, we have

\begin{cor}\label{convergence.uk.2}
Let $\{u_k\}$ satisfy the same conditions as in Lemma 3.3. Then, for any fixed $\theta$, $u_k(s,\theta)$ converges to $\gamma$
in $C^1[0,a]$. Moreover, we have
$$t_k^a\rightarrow\infty,\s\s \sqrt{|\mbox{Re} \ \beta(u_{k})|}t_k^a <C(a),$$
and
\begin{equation}\label{arclength2}
\left\|\frac{1}{\sqrt{| \mbox{Re}(\beta(u_{k}))|}}\left|\frac{\partial s}{\partial t}(s)\right|-\frac{1}{\sqrt{2\pi}}\right\|_{C^0([0,a])}\rightarrow 0.
\end{equation}
Here $C(a)$ is a positive constant which depends on $a$.
\end{cor}

\medskip
Now we turn to the discussions on the asymptotic behavior of the index and the second variation of the energy of $u_k$. Since $\mbox{Ric}_N\geq\lambda>0$, by the well-known Myers theorem we know that, if
 $$a\geq\frac{\pi}{\sqrt{\lambda(n-1)^{-1}}}+2\epsilon,$$
then there exists a tangent vector field $V_0(s)$ on $N$, which is smooth on $\gamma$, and 0 on $\gamma|_{[0,\epsilon]}$ and $\gamma|_{[a-\epsilon,a]}$,
such that the second variational of length of $\gamma$ satisfies
\begin{equation}
I_\gamma(V_0,V_0)=\int_0^a(\langle \nabla_{\dot{\gamma}}V_0,\nabla_{\dot{\gamma}}V_0 \rangle-R(V_0,\dot{\gamma},\dot{\gamma},V_0))ds<-\delta<0.
\end{equation}

Following the arguments in Section 2, we can see easily that the conclusions in Proposition \ref{sy}
are implied by the following Lemma.
\begin{lem}
Let $\{u_k\}$ satisfy the same conditions as in Lemma 3.3. Then, for sufficiently large $k$, there exists a section $V_k$ of
$u_{k}^*(TN)$, which is supported in $[0,t_k^a]$, such that
$$\delta^2E(u_k)(V_k, V_k)<0.$$
\end{lem}

\proof Let $P$ be projection from $T\R^K$ to $TN$. We define
 $$V_k(t,\theta)=P_{u_k(s,\theta)}(V_0(s)),$$
where $s$ is the arc-length parametrization of $u_{k}^*(t)$ with $s(0)=0$. Then, $V_k$ is a smooth section of $u_{k}^{-1}(TN)$ which
is supported in $[0,\, t_k^a]\times S^1$. By Corollary \ref{convergence.uk.2}, for any fixed $\theta$, we have
$$V_k(u_k(se^{\sqrt{-1}\theta}))\rightarrow V_0(s)\s\s \mbox{in}\s C^1[0,\,a].$$

Next, we will show
\begin{equation}\label{I.I2}
\lim_{k\rightarrow+\infty}
\frac{1}{\sqrt{|\mbox{Re} \ \beta(u_{k})|}}\delta^2E(u_k)(V_k,V_k)
=2\sqrt{2\pi}I_\gamma(V_0,V_0).
\end{equation}
Since
\begin{eqnarray*}
\delta^2E(u_k)(V_k,V_k)&=&2\int_0^{2\pi}\int_{0}^{t_k^a}
(\langle \nabla V_k,\nabla V_k \rangle-R(V_k,\nabla u_k,\nabla u_k,V_k))dtd\theta\\
\\
&=&
2\int_0^{2\pi}\int_{0}^{t_k^a}
(\langle \nabla_{\frac{\partial u_k}{\partial t}} V_k,\nabla_{\frac{\partial u_k}{\partial t}} V_k \rangle-R(V_k,\frac{\partial u_k}{\partial t},\frac{\partial u_k}{\partial t},V_k))dtd\theta\\
&&+2\int_0^{2\pi}\int_{0}^{t_k^a}
(\langle \nabla_{ \frac{\partial u_k}{\partial \theta}} V_k,\nabla_{\frac{\partial u_k}{\partial \theta}} V_k \rangle-R(V_k, \frac{\partial u_k}{\partial \theta},\frac{\partial u_k}{\partial \theta},V_k))dtd\theta\\
&=&
2\mathbf{I}+2\mathbf{II}.
\end{eqnarray*}
Noting
\begin{eqnarray*}
\frac{\mathbf{I}}{\sqrt{|\mbox{Re} \ \beta(u_{k})|}}&=&
\int_0^{2\pi}\int_{0}^{a}\left(\langle \nabla_{\frac{\partial u_k}{\partial s}} V_k,\nabla_{\frac{\partial u_k}{\partial s}} V_k \rangle-R(V_k,\frac{\partial u_k}{\partial s},\frac{\partial u_k}{\partial s},V_k)\right)\\
&&\times\frac{\left|\frac{\partial s}{\partial t}\right|}{\sqrt{|\mbox{Re} \ \beta(u_{k})|}}dsd\theta,
\end{eqnarray*}
we infer from Corollary \ref{convergence.uk.2}
$$\lim_{k\rightarrow+\infty}\frac{\mathbf{I}}{\sqrt{|\mbox{Re} \ \beta(u_{k})|}}
=\sqrt{2\pi}I_\gamma(V_0,V_0).$$

On the other hand, we have
\begin{eqnarray*}
\frac{\mathbf{II}}{\sqrt{|\mbox{Re} \ \beta(u_{k})|}}&=&
\frac{1}{\sqrt{|\mbox{Re} \ \beta(u_{k})|}}\int_0^{2\pi}\int_{0}^{t_k^a}
\left(\langle \nabla_{ \frac{\partial u_k}{\partial \theta}} V_k,\nabla_{\frac{\partial u_k}{\partial \theta}} V_k \rangle-R(V_k,\frac{\partial u_k}{\partial \theta},\frac{\partial u_k}{\partial \theta},V_k)\right)dtd\theta\\
&\leq&
\frac{C}{\sqrt{|\mbox{Re} \ \beta(u_{k})|}}\int_0^{2\pi}\int_{0}
^{t_k^a}|\frac{\partial u_k}{\partial \theta}|^2dtd\theta.
\end{eqnarray*}
For any given $T>0$, we set $$m_k=\left[\frac{t_k^a}{T}\right]+1.$$
By Corollary \ref{convergence.uk.2}, there holds
$$\sqrt{|\mbox{Re} \ \beta(u_{k})|}m_k\leq C(T).$$
Hence, it follows
\begin{eqnarray*}
\frac{\mathbf{II}}{\sqrt{|\mbox{Re} \ \beta(u_{k})|}}
&\leq&
\frac{C}{\sqrt{|\mbox{Re} \ \beta(u_{k})|}}\int_{\cup_{i=0}^{m_k}[iT,(i+1)T]\times S^1}
|u_{k,\theta}|^2dtd\theta\\
&\leq&
\frac{Cm_k\sqrt{|\mbox{Re} \ \beta(u_{k})|}}{ |\mbox{Re} \ \beta(u_{k})| }\frac{1}{m_k}\int_{\cup_{i=0}^{m_k}[iT,(i+1)T]\times S^1}
|u_{k,\theta}|^2dtd\theta\\
&\leq&
\frac{C(T)}{m_k}\left(\frac{1}{|\mbox{Re} \ \beta(u_{k})|}\int_{\cup_{i=0}^{m_k}[iT,(i+1)T]\times S^1}
|u_{k,\theta}|^2dtd\theta\right).
\end{eqnarray*}
In view of Lemma \ref{theta.energy2}, we concludes
$$
\lim_{k\rightarrow\infty}\frac{1}{\sqrt{|\mbox{Re} \ \beta(u_{k})|}}\mathbf{II}=0.
$$
Immediately, it follows
$$\lim_{k\rightarrow+\infty}
\frac{1}{\sqrt{|\mbox{Re} \ \beta(u_{k})|}}\delta^2E(u_k)(V_k,V_k)
=2\sqrt{2\pi}I_\gamma(V_0,V_0).$$
Hence, for $k$ large enough, we have the desired inequality
$$\delta^2E(u_k)(V_k, V_k)<0.$$
Thus, we complete the proof of this lemma.
\endproof

\medskip

\noindent{\bf Acknowledgement:} Y. Li supported by NSFC (Grant No. 11131007), Y. Wang supported by NSFC(Grant No. 11471316).


{}

\vspace{1.0cm}

Yuxiang Li

{\small\it Department of Mathematical Sciences, Tsinghua University,  Beijing 100084, P.R.China.}

{\small\it Email: yxli@math.tsinghua.edu.cn.}\\

Lei Liu

{\small\it Department of Mathematical Sciences, Tsinghua University, Beijing 100084, P.R.China.}

{\small\it Email: liulei1988@mail.tsinghua.edu.cn }\\

Youde Wang

{\small\it Academy of Mathematics and Systems Sciences, Chinese Academy of Sciences, Beijing 100080,  P.R. China.}

{\small\it Email: wyd@math.ac.cn}

\end{document}